\documentclass[11pt, a4paper]{amsart} 
\usepackage[utf8]{inputenc}
\usepackage[T1]{fontenc}
\usepackage{listings}
\usepackage{xcolor}
\usepackage{tikz}
\usepackage{pgfplots}
\usepackage{adjustbox}
\usepackage{array}
\usepackage{caption}
\usetikzlibrary{er,positioning}
\usetikzlibrary{calc,intersections,through,backgrounds}

\usepackage{pdflscape}

\definecolor{codegreen}{rgb}{0,0.6,0}
\definecolor{codegray}{rgb}{0.5,0.5,0.5}
\definecolor{codepurple}{rgb}{0.58,0,0.82}
\definecolor{backcolour}{rgb}{0.95,0.95,0.92}
 
\lstdefinestyle{mystyle}{
    backgroundcolor=\color{backcolour},   
    commentstyle=\color{codegreen},
    keywordstyle=\color{magenta},
    numberstyle=\tiny\color{codegray},
    stringstyle=\color{codepurple},
    basicstyle= \tiny,  
    breakatwhitespace=false,         
    breaklines=true,                 
    captionpos=b,                    
    keepspaces=true,                 
    numbers=left,                    
    numbersep=5pt,                  
    showspaces=false,                
    showstringspaces=false,
    showtabs=false,                  
    tabsize=1
}
 
\usepackage[utf8]{inputenc}
\usepackage[T1]{fontenc}
\usepackage{listings}
\usepackage{amssymb,epsfig}
\usepackage{verbatim}
\usepackage{fancyvrb}

\usepackage[T1]{fontenc}
\usepackage{epsfig,multicol, graphics}
\usepackage{hyperref}
\usepackage{pdfpages}
\usepackage{marginnote}
\usepackage[all]{xy}
\usepackage{graphics}
\usepackage{graphicx}
\usepackage{color}

\newtheorem{thm}{Theorem}
\newtheorem{pro}[thm]{Proposition}

\newtheorem{lem}[thm]{Lemma}

\newtheorem{exe}[thm]{Example}

\newcommand{\E}{{\mathcal E}}
\newcommand{\F}{{\mathcal F}}

\newcommand{\x}{{\bf x}}

\newcommand{\plano}{{\mathbb P}^2}

\newcommand{\tor}{\xymatrix{\ar@{-->}[r]&}}
 
\definecolor{grass}{rgb}{0.14,0.72,0.2}

\begin{document}

\lstset{language=Python}

{\tiny \title{Negative Curves and Elliptic Fibrations on a Special Rational Surface}}
\author{by\,   Lu\'is Gustavo Mendes and Liliana Puchuri}

\begin{abstract}

The blown up complex projective plane in the twelve triple points of the dual Hesse arrangement  has an  infinite number of irreducible rational curves of self-intersection $-1$, for short, $(-1)$-curves. In the preprint version of  \cite{Dumnicki}, T. Szemberg \emph{et alii}  tried to  keep  track of   plane rational curves which produce  $(-1)$-curves,  by applying successively to a straight line compositions taken from   three different quadratic Cremona maps preserving the dual Hesse  arrangement, in such a way  as to produce  a  diagram of bifurcations. The symmetric aspect of the degrees of curves encoded in the entries of the   diagram motivated them to  propose the  problem    of \emph{to give a closed formula for the degree of the rational curve at the entry $a_{i,j}$}. 
We solved  this  problem for a slightly different diagram of bifurcations of the same Cremona maps, which recover and  extends the data on degrees of the original diagram. We also are able to describe the position and multiplicites of singularities of each rational curve encoded in the diagram, as well as their equations (implemented in \emph{Python} and \emph{Singular} softwares).   
Our   idea is  to embed the  plane rational curves as components  of special elements of elliptic pencils, then apply a closed  formula      for the degrees of the generic curves of the elliptic pencils of  \cite{Puchuri}, and,  from this,  extract   the data  of the  rational curves on each entry of the diagram.

\end{abstract}

\maketitle

\section{Introduction and results}

In the complex projective plane there are no projective curves with negative self-intersection, thanks to Bézout Theorem. Each time we blow up the plane we introduce at least one  $(-1)$-curve, namely, the exceptional divisor of the blow up. When we blow up the plane in two distinct points, besides the two exceptional divisors, the  projective line connecting the points becomes a $(-1)$-curve of the rational surface obtained, denoted $Bl_2 (\plano)$.

A non trivial fact  is that, after  blow ups of \emph{nine or more} distinct points in the plane, the rational surface  $Bl_{n\geq 9}(\plano)$    may  have an infinite number of $(-1)$-curves.
In \cite{Dumnicki}, T. Szemberg et alli considered negative curves on a special rational surface wich arises after blow ups of nine among  the twelve   points belonging to   the \emph{dual Hesse arrangement} (see Section \ref{hesse}). 

In a  preprint version  \cite{preprint}, they   proposed to track rational plane curves which give rise to  $(-1)$-curves on the special rational surface.
They started  with a straight line and applied to it one among three options of quadratic Cremona mas $\phi_1, \phi_2, \phi_3$,  where each $\phi_i$ preserves the  dual Hesse arrangement (in the sense of  strict transform). Being involutive, it has no use to apply two times the same $\phi_i$, but  at each stage there are two options  to be applied   on  the curve. Of course, different paths can achive at the same results.  In this way they produced a    diagram of bifurcations of Cremona maps  and their effect on degrees of rational curves.      
Problem 3.10 of \cite{preprint} includes \emph{to find a direct formula for the entry in line $i$ column $j$  }.

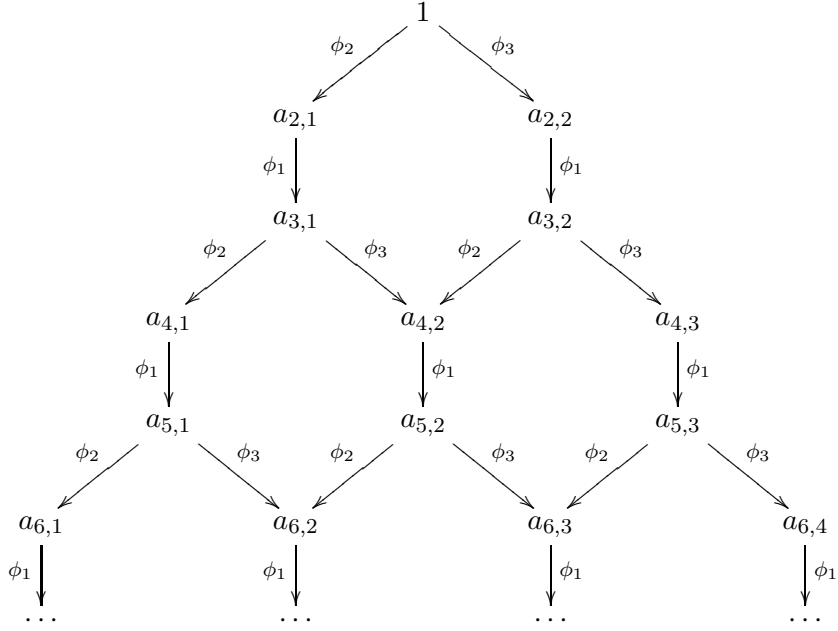
\begin{figure}
	\centerline{
		\xymatrix{
			&  &  & {1} \ar[ld]_{\phi_2}\ar[rd]^{\phi_3} &  &  &  \\
			&  & a_{2,1}\ar[d]_{\phi_1} &  & a_{2,2}\ar[d]^{\phi_1} &  &  \\
			&  & a_{3,1} \ar[ld]_{\phi_2}\ar[rd]^{\phi_3}&  & a_{3,2}\ar[ld]_{\phi_2}\ar[rd]^{\phi_3} &  &  \\
			& a_{4,1}\ar[d]_{\phi_1} &  & a_{4,2}\ar[d]^{\phi_1} &  & a_{4,3}\ar[d]^{\phi_1} &  \\
			& a_{5,1}\ar[ld]_{\phi_2}\ar[rd]^{\phi_3} &  & a_{5,2}\ar[ld]_{\phi_2}\ar[rd]^{\phi_3} &  & a_{5,3}\ar[ld]_{\phi_2}\ar[rd]^{\phi_3} &  \\
			a_{6,1}\ar[d]_{\phi_1} &  & a_{6,2}\ar[d]_{\phi_1} &  & a_{6,3}\ar[d]^{\phi_1} &  & a_{6,4}\ar[d]^{\phi_1} \\
			{\cdots} &  & {\cdots} &  & {\cdots} &  & {\cdots}
	}}
\caption{First stages of    a sequence of bifurcations of quadratic Cremona maps.}
\label{diagramours}
\end{figure}

Our  ideas for solving this problem  are the following: 

\begin{itemize}

\item i) embed the rational curves   as components of  special elements of elliptic pencils $\E_t$; 
\item ii) apply to the elliptic pencils  quadratic Cremona maps preserving the dual Hesse arrangement, forming a  diagram of bifurcations (Figure  \ref{diagramours})  slightly different from the original diagram,  but which   recovers and extends the data  of the original one;
\item iii) at each stage, the  rational curve becomes a $(-1)$-curve after part  of the  blow ups at  twelve fixed  special points; the blow up of the twelve points produces  from $\E_t$  a non-minimal elliptic fibration,   containing the  $(-1)$-curve as  a  component of a singular  fiber;  
\item iv) after each application of $\phi_i$, we control the change in the parameter $t$ of $\E_t$ and determine  $t$ in  terms of the entry $a_{i,j}$  of  the diagram;  
\item  v) use the  formula of \cite{Puchuri} for the degree   of the generic element of the elliptic pencil  $\mathcal{E}_{t}$;
\item  vi) extract data of   degree, singularities and even the equation of   the rational curves $Cr_{i,j}$ corresponding to the entry  $a_{i,j}$.
\end{itemize}

Our  result is:
\begin{thm}\label{fatoprincipal}
Consider a diagram of any  number of  bifurcations of the three quadratic Cremona maps $\phi_ 1,\phi_2,\phi_3$  following repeatedly the scheme of bifurcations of  Figure    \ref{diagramours}.  The $\phi_i$ are the same as   of T. Szemberg et alli (described  in  Section \ref{hesse}). The diagram encodes the degrees of the  strict transforms of rational plane curves.  It begins with the degree $1$ of  a straight line; in the  second and third  rows it has two columns; for all $i \geq 4$, it has,  in the $i$-ith row,   $\frac{i}{2} +1$ columns, if $i$ is even, or  $\frac{i-1}{2} +1$ columns, if $i$ is odd. 

Let $Cr_{i,j}$ be a rational curve corresponding  to position $a_{i,j}$ in such diagram. Then 
 for  all  $i  \geq 2$ 
\[\mbox{deg}(Cr_{i,j})= c^2(i,j)+ d^2(i,j) + c(i,j) \cdot d(i,j)- c(i,j) + d(i,j)
\]
where

\begin{itemize}

\item  for even $i\geq 2$  and  $1  \leq    j \leq  \frac{i}{2} +1$, 
\[c(i,j)= 1-j, \quad d(i,j) = \frac{i}{2} + 1- j \]
\item for odd   $i\geq 2$  and  for $1  \leq    j \leq \frac{i-1}{2}  +1$,  
\[  c(i,j)= j,\quad   d(i,j)= j -\frac{i-1}{2} -2  \]
\end{itemize}
\end{thm}

The first stages   produce a sequence of degrees recovering and extending the degrees of the original diagram:
\begin{center}
\begin{tabular}{>{$}l<{$\hspace{12pt}}*{13}{c}}
&&&&&&&\color{blue} 1&&&&&&\\
&&&&&&\color{blue} 2&&\color{blue} 2&&&&&\\
&&&&&&\color{blue} 4&&\color{blue} 4&&&&&\\
&&&&& \color{blue} 6 &&\color{blue} 5&& \color{blue} 6 &&&&\\
&&&&& \color{blue} 9 &&\color{blue} 8&& \color{blue} 9 &&&&\\
&&&& \color{blue} 12 &&\color{blue} 10&& \color{blue} 10 && \color{blue} 12  &&&\\
&&&& \color{blue} 16 &&\color{blue} 14&& \color{blue} 14 && \color{blue} 16  &&&\\
&&& \color{blue} 20 &&\color{blue} 17&& \color{blue} 16 && \color{blue} 17 &&\color{blue} 20 &&\\
&&& \color{blue} 25 &&\color{blue} 22&& \color{blue} 21 && \color{blue} 22 &&\color{blue} 25 &&\\
&& \color{blue} 30 &&\color{blue} 26&& \color{blue} 24 && \color{blue} 24 &&\color{blue} 26  &&\color{blue} 30 &\\
&& \color{blue} 36 &&\color{blue} 32&& \color{blue} 30 && \color{blue} 30 &&\color{blue} 32  &&\color{blue} 36 &\\
\end{tabular}
\end{center}
\centerline{\color{blue} etc.}

As a second result,  we  give the positions in the plane (explained in Section \ref{hesse})  and multiplicities of the  singular points of the rational curves  $Cr_{i,j}$  at any entry of the diagram in  Figure \ref{diagramours}. 

\newpage

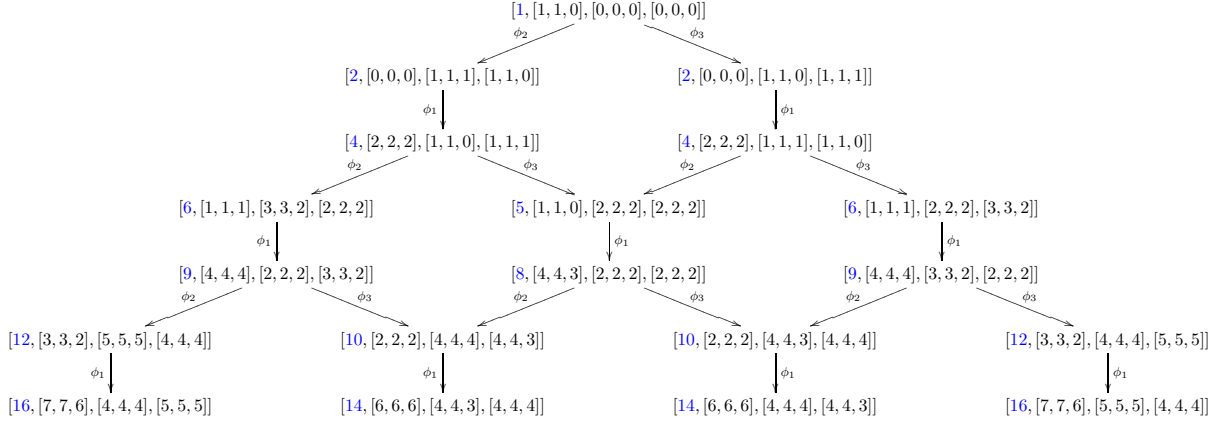
\begin{figure}
	\centerline{
		\scalebox{0.6}{
			\xymatrix@C=-1cm{
			&  &  & [\textcolor{blue}{1},[1,1,0],[0,0,0],[0,0,0]] \ar[ld]_{\phi_2}\ar[rd]^{\phi_3} &  &  &  \\
			&  & [\textcolor{blue}{2},[0,0,0],[1,1,1],[1,1,0]]\ar[d]_{\phi_1} &  & [\textcolor{blue}{2},[0,0,0],[1,1,0],[1,1,1]]\ar[d]^{\phi_1} &  &  \\
			&  & [\textcolor{blue}{4},[2,2,2],[1,1,0],[1,1,1]] \ar[ld]_{\phi_2}\ar[rd]^{\phi_3}&  & [\textcolor{blue}{4},[2,2,2],[1,1,1],[1,1,0]]\ar[ld]_{\phi_2}\ar[rd]^{\phi_3} &  &  \\
			& [\textcolor{blue}{6},[1,1,1],[3,3,2],[2,2,2]]\ar[d]_{\phi_1} &  & [\textcolor{blue}{5},[1,1,0],[2,2,2],[2,2,2]]\ar[d]^{\phi_1} &  & [\textcolor{blue}{6},[1,1,1],[2,2,2],[3,3,2]]\ar[d]^{\phi_1} &  \\
			& [\textcolor{blue}{9},[4,4,4],[2,2,2],[3,3,2]]\ar[ld]_{\phi_2}\ar[rd]^{\phi_3} &  & [\textcolor{blue}{8},[4,4,3],[2,2,2],[2,2,2]]\ar[ld]_{\phi_2}\ar[rd]^{\phi_3} &  & [\textcolor{blue}{9},[4,4,4],[3,3,2],[2,2,2]]\ar[ld]_{\phi_2}\ar[rd]^{\phi_3} &  \\
			[\textcolor{blue}{12},[3,3,2],[5,5,5],[4,4,4]]\ar[d]_{\phi_1} &  & [\textcolor{blue}{10},[2,2,2],[4,4,4],[4,4,3]]\ar[d]_{\phi_1} &  & [\textcolor{blue}{10},[2,2,2],[4,4,3],[4,4,4]]\ar[d]^{\phi_1} &  & [\textcolor{blue}{12},[3,3,2],[4,4,4],[5,5,5]]\ar[d]^{\phi_1} \\
[\textcolor{blue}{16},[7,7,6],[4,4,4],[5,5,5]] &  & [\textcolor{blue}{14},[6,6,6],[4,4,3],[4,4,4]] &  & [\textcolor{blue}{14},[6,6,6],[4,4,4],[4,4,3]] &  & [\textcolor{blue}{16},[7,7,6],[5,5,5],[4,4,4]]
	}}
}
\caption{Bifurcations of quadratic Cremona maps $Q_i$ applied to irreducible  rational curves. Degrees in blue and data of singularities.}
\label{diagramrational1}
\end{figure}

For example,   $Cr_{9, 3}$ is a rational curve of degree $21$ with ordinary singular  points,  namely,   two  $9$-uple  points, one $8$-uple point  and six $6$-uple points whose positions we can describe.

Even the equations of the rational curves in each entry are available,  with  a code in \emph{Singular} at the end of  this paper. For instance,
{\small 
\[Cr_{3,1}:  x^4 -(\tau+1) x^3y-(2\tau +1) x^2y^2- (\tau +1) xy^3+
  y^4+\tau x^3z+(\tau -1) x^2yz- (\tau +1) xy^2z+\]
\[+ \tau  y^3z+(\tau +2) x^2z^2 +(\tau +2) xyz^2+(\tau +2) y^2z^2-
(\tau + 1) xz^3-(\tau +1) yz^3-(\tau +1) z^4 = 0\]}
where $\tau $ is the primitive 3th root of unit.

\tableofcontents

\pagebreak 
\section{Dual Hesse arrangement and three quadratic Cremona maps preserving it}\label{hesse}

The  \emph{dual Hesse arrangement} of projective lines and points   on the complex projective plane 
is composed by a set of nine    lines $\mathcal{L}_9$  intersecting at  twelve  points $\mathcal{P}_{12}$. Each one of the twelve points  is a triple point of  the set of  lines  and  on  each line there are four of the twelve points. In the research field of arrangements these incidences are   usualy denoted  by $(12_3, 9_4  )$.

The   lines of $\mathcal{L}_9$  can be  given  in coordinates $(x:y:z)$ by 
\begin{equation}\label{lines}
\begin{cases} 
l_1:= y-x =0;\quad  l_2:= y-\tau \cdot  x =0;\quad   l_3:= y-\tau^2 \cdot x =0 \\ 
  m_1:= z-x =0;\,\,  m_2:= z- \tau \cdot x =0; \,\, m_3:= z-\tau^2 \cdot x=0\\
 n_1:= z-y=0;\quad  n_2:= z-\tau \cdot y =0;\,\,  n_3:= z-\tau^2 \cdot y =0;
\end{cases}
\end{equation}
where $\tau$ is a primitive 3th root of unit; an  equation for all  the set $\mathcal{L}_9$ is 
\[(x^3 - z^3) \cdot (y^3-z^3) \cdot (x^3-y^3) = 0 \]
 
The set of twelve points $\mathcal{P}_{12}$ can be decomposed  in sets of three points 
\[\mathcal{P}_{12} = \mathcal{P}_3 (1)\cup \mathcal{P}_3 (\tau) \cup \mathcal{P}_{3}(\tau^2) \cup \mathcal{P}_{3} (\infty)\]  
where 
\[
(\star) \begin{cases} 
\mathcal{P}_3(1):\,  (1:1:1)= l_1 \cap m_1\cap n_1,\quad  (1:\tau:\tau^2)= l_2\cap m_3\cap n_2,\quad  (1:\tau^2:\tau) = l_3\cap m_2\cap n_3 \\ 
\mathcal{P}_3 (\tau):\,  (1:\tau:1)= l_2\cap m_1\cap n_3 ,\quad  (1:1:\tau)= l_1\cap m_2\cap n_2 ,\quad    (\tau:1:1)=  l_3 \cap m_3 \cap n_1   \\
\mathcal{P}_3 (\tau^2):\, (1:1:\tau^2)= l_1\cap m_3 \cap n_3,\,\,  (\tau^2:1: 1) =  l_2 \cap m_2 \cap n_1 \,\,\, \,  \,   (1:\tau^2:1) = l_3 \cap m_1 \cap n_2 \\
\mathcal{P}_3 (\infty):\,  (0:0:1) = \, l_1 \cap l_2 \cap l_3,\quad \,   (0:1:0)= m_1 \cap m_2 \cap m_3   ,\,\,    (1:0:0)= n_1 \cap n_2 \cap n_3  
\end{cases}\]

The dual  Hesse arrangement is projectively rigid (cf.  \cite{LinsNeto} Prop. 1 or  \cite{Magdalena} Th. 1): that is,  the incidences  $(12_3, 9_4  )$ determine the dual Hesse arrangement up to automorphism of the complex projective plane.

\begin{figure}
\begin{tikzpicture}[scale=1.5]
	\coordinate (O) at (0,0);
	\draw[thick, name path=O15] (O)-- (15:5);	
	\draw[thick, name path=O30] (O) -- (30:5);
	\draw[thick, name path=O45] (O) -- (45:5);
	\coordinate (V1) at (5,0);	
	\draw[thick, name path=V115] (V1) -- +(135:5);
	\draw[thick, name path=V130] (V1) -- +(150:5);
	\draw[thick, name path=V145] (V1) -- +(165:5);
	\coordinate (V2) at (60:5cm);
	\draw[thick,blue, name path=V215] (V2) -- +(255:5);
	\draw[thick, name path=V230] (V2) -- +(270:5);
	\draw[thick,green!50!black, name path=V245] (V2) -- +(285:5);
	\path [name intersections={of=O15 and V115,by=A1}];
	\path [name intersections={of=O30 and V115,by=A2}];
	\path [name intersections={of=O45 and V115,by=A3}];
	\path [name intersections={of=O15 and V215,by=B1}];
	\path [name intersections={of=O30 and V215,by=B2}];
	\path [name intersections={of=O45 and V215,by=B3}];
	\path [name intersections={of=O30 and V130,by=C1}];
	\path [name intersections={of=O15 and V145,by=C2}];
	\path [name intersections={of=O15 and V130,by=C3}];
	\path [name intersections={of=V215 and V115,by=D1}];
	\path [name intersections={of=V245 and V145,by=D2}];
	\path [name intersections={of=V245 and O45,by=D3}];
	\path [name intersections={of=O45 and V145,by=D4}];
	\draw[blue, thick] ([shift={(-105:0.3cm)}]A1.center) -- ([shift={(-105:-0.5cm)}]A1.center);
	\draw[green!50!black, thick] ([shift={(-75:0.3cm)}]D4.center) -- ([shift={(-75:-0.5cm)}]D4.center);
	\draw[orange,fill] (O) circle (2.5pt);
	\draw[orange,fill] (V1) circle (2.5pt);
	\draw[orange,fill] (V2) circle (2.5pt);
	\draw[orange,fill] (A1) circle (2.5pt);
	\draw[orange,fill] (A2) circle (2.5pt);
	\draw[orange,fill] (A3) circle (2.5pt);
	\draw[white,fill] (B1) circle (1.5pt);
	\draw[orange,fill] (B2) circle (2.5pt);
	\draw[orange,fill] (B3) circle (2.5pt);
	\draw[orange,fill] (C1) circle (2.5pt);
	\draw[orange,fill] (C2) circle (2.5pt);
	\draw[orange,fill] (C3) circle (2.5pt);
	\draw[white,fill] (D1) circle (1.5pt);
	\draw[white,fill] (D2) circle (1.5pt);
	\draw[white,fill] (D3) circle (1.5pt);
	\draw[orange,fill] (D4) circle (2.5pt);
\end{tikzpicture}
\caption{In this real figure of the dual Hesse arrangement  two  lines (blue and green)  are poorly represented as broken and nonconnect. The are no double intersections, just triple intersections of lines at orange points.}
\label{Dual}
\end{figure}
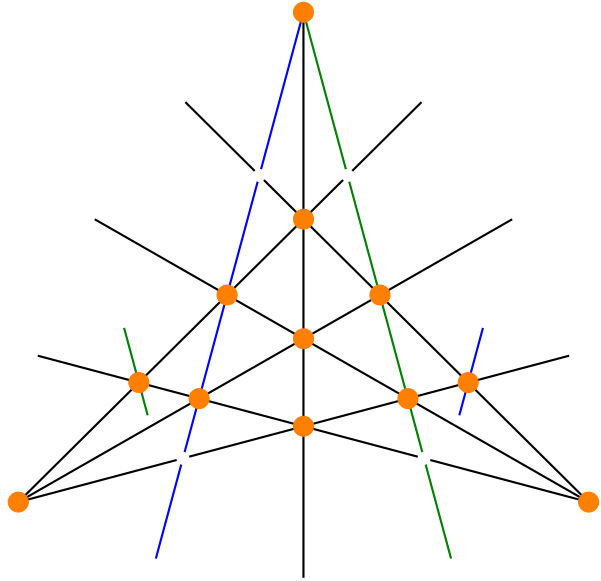

Any attempt to illustrate  the dual Hesse arrangement on the \emph{real} plane  has some deficiency: either some  lines shall be represented as curved, or as being  not-connected, or  some points of the arrangement will  be missing. We  illustrate it in Figure \ref{Dual}.

There are  projective automorphisms that preserve the dual Hesse arrangement, but there also quadratic Cremona maps preserving it, in the sense that the strict transform of the set of lines is preserved as well  as the incidences $(12_3, 9_4  )$.

We shall use the following three:  $Q_1$  with  indetermination points at $\mathcal{P}_3 (1)$,  $Q_\tau$ with indetermination points at $\mathcal{P}_3 (\tau)$ and $Q_{\tau^2}$  with indetermination points at $\mathcal{P}_{\tau^2}$; in particular, no line of the fundamental triangles of these three quadratic maps belongs to $\mathcal{L}_9$.

In homogeneus coordinates $(x:y:z)$: 
\[(\star\star)\begin{cases}   Q_1(x:y:z) = (y^2 -x z : x^2 - y z: z^2- x y  ) \\ Q_\tau (x:y:z)= (\tau y^2 - x z, \tau x^2 - y z: z^2- \tau^2 x y    ) \\ Q_{\tau^2} (x:y:z)= (\tau^2 y^2 - x z, \tau^2 x^2 - yz : z^2- \tau x y)  \end{cases}   \]

In   notations of \cite{preprint},    \[Q_1 = \phi_1, \quad  Q_\tau =  \phi_2,\quad Q_{\tau^2}= \phi_3\]

\section{Elliptic pencils containing $\mathcal{L}_9$ as components of special elements}

Consider a pencil  of curves  in the complex projective plane $\plano$  whose generic element is irreducible and has geometrical genus $g$; for $g=1$, let us represent  such an  \emph{elliptic pencil} in the plane as 
\[\E:\quad  c_1 \cdot f(x:y:z) + c_2 \cdot g(x:y:z) = 0, \quad (c_1:c_2) \in \mathbb{P}^1,\]
where $f, g$ are homogeneous polynomials of the same degree. 

We can   associate to it a \emph{singular holomorphic foliation} of $\plano$, with finite singular set, as follows. From  the global rational function $\Psi: \plano \dashrightarrow \mathbb{P}^1$  given by
\[\Psi =  \frac{f}{g} = \frac{-c_2}{c_1}=:c   \]
we produce a differential equation
\[0 = d ( \frac{f}{g} ) =  \frac{g df - f dg}{g^2}\]
We  start by elliminating  the curve of poles, getting 
\[g \, df - f \,  dg = 0\]
But the left side of this equation may have curves of zeroes, which must be elliminated also. In fact, let  
\[ C_s = {\prod_{s, k}} (C_{s,k})^{\alpha_ {s,k}}  \]
 be a decomposition in irreducible factors of special elements $C_s$ of $\E$ (in finite number thanks to Bertini's theorem). Then the zero divisor of $ g \, df - f \,  dg $ 
has degree   \[\sum_{s, k}  (\alpha_{s,k} -1)\cdot  \mbox{deg}(C_s)\] 
After elliminating divisors  of zeroes and poles of $ g \, df - f \,  dg $   we get a foliation with finite singular set:
\[\Gamma =  a(x:y:z) \, dx + b(x:y:z) \, dy + c(x:y:z) \, dz = 0\] 
A fundamental remark (which goes back to G. Darboux, H. Poincaré, P. Painlevé) is  that, although the degree of the generators  $f$ and $g$ of $\E$ can be high, the degree of the polynomial coeficients $a,b,c$ of $\Gamma$ can be low. This happens when  there are  multiple components of special elements  (i.e. with  $\alpha_{s,k }\geq 2$) of high multiplicity or having  high degree of is support. 

In the inverse direction we have non-trivial  question: given a  foliation with finite singular set  of  the plane, to  determine if it corresponds to a pencil of curves (i.e. if it  has a global first integral)  and, in the positi~ve case, to determine the geometric genus of the generic element of the pencil.

\section{Lins Neto's foliations and rational components of elliptic pencils}\label{rationalelliptic}

In \cite{LinsNeto}, A. Lins Neto produced  a $1$-parameter family of singular holomorphic foliations with finite singular set
\[\Gamma_t = \Omega + t \cdot \Xi  = 0,\quad t\in \mathbb{P}^1\]
where 
$\Omega$ and $\Xi $ are  polynomial $1$-forms whose coefficients have  degree five. with the following remarkable properties:

\begin{itemize}

\item For all  $t \in \mathcal{P}^1$, Lins Neto's foliation $\F_t$ is tangent to the set of lines  $\mathcal{L}_9$ of the dual Hesse arrangement.

\item Just for  a countable and dense set of the parameters  $D \subset \mathbb{P}^1$, the foliation   $\Gamma_{t \in D} $ corresponds to an elliptic pencil $\E_{t\in D}$.

\item For all  $t \in \mathcal{P}^1 \setminus \{1,\tau,\tau^2, \infty\}$,  the singularities of $\F_t$ are composed by twelve $(1,1)$-type   points located at the set  $\mathcal{P}_{12}$ of the dual Hesse arrangement and nine $(1,-3)$-type mobile singularites $p_1(t), \ldots p_{9}(t)$,  belonging to $\mathcal{L}_9 \setminus \mathcal{P}_{12}$, one at each line of $\mathcal{L}_9$. By $(1,1)$-type  or \emph{radial} point we mean a singularitity of foliation with  \emph{local} meromorphic first integral $\frac{y}{x} = c$;  by $(1,-3)$-type we mean  singularity of foliation with \emph{local} analytic first integral $y\cdot x^3 = c$.

\item After the  work of Lins Neto and  M. McQuillan, A. Guillot and  L. Puchuri (\cite{Mc},  \cite{Guillot}, \cite{Puchuri}) it is  known that for all 
$t \in \mathbb{Q}(\tau)$ the foliation $\F_{p+q \tau}$ corresponds to an elliptic pencil, denoted $\E_{p + q \tau}$. 

\item If $\pi(n)$ counts the number of elements in $\mathbb{Q}(\tau)$ for which $\E_t$ have  generic curves of degree  at most $n$, then $\pi(n) = O(n^2)$ (Corollary 9 of \cite{Puchuri})

\end{itemize}

All  elliptic pencils $\E_{p+q\tau}$ of lins Neto's family  have special elements containing lines of $\mathcal{L}_9$.  But they have also components of special elements which are irreducible rational curves with multiplicity $\alpha_{s,k} =3$.

A result of  \cite{MPEffective} is that,    for any $t = m+n\tau \in \mathbb{Z}(\tau)$,  the elliptic pencil $\E_{m+n \tau}$ can be represented as 
\[\E_{m+n\tau}: \quad   c_1 \cdot P_3 \cdot  (Cr_1)^3 + c_2 \cdot  Q_3 \cdot  (Cr_2)^3  = 0\]
where 

\begin{itemize}

\item $P_3=0$ and  $Q_3=0$ define sets of three distinct straight lines belonging to $\mathcal{L}_9$; 

\item  $Cr_1 =0$ and $Cr_2=0$  define two irreducible rational curves of the same degree;

\item there is  a third special element  $R_3 \cdot (Cr_3)^3 $, for  $R_3=0$  a set of three lines of $\mathcal{L}_9$ and $Cr_3$ an irreducible rational curve.

\item $P_3 \cdot Q_3 \cdot R_3 = 0$ coincides with $\mathcal{L}_{12}$ 

\end{itemize}

Each  elliptic pencil $\E_{m+n\tau}$  produces, after blow up of the twelve points $\mathcal{P}_{12}$ of the dual Hesse arrangement, a non-minimal elliptic fibration of the blown up plane in twelve points $Bl_{12}(\plano)$:
\[F_{m+n\tau} : Bl_{12}(\plano) \rightarrow \mathbb{P}^1\]
Each fibration $ F_{m+n\tau}$ has three special fibres $\tilde{IV}$, isomorphic to  the Kodaira type $IV$ blown up at its  triple point. That is, each special fiber has a "central" component which is a  $(-1)$-curve with  multiplicity three and three components with self intersection $(-3)$ and multiplicity one.

The three  "central"  components are  the strict transforms in $Bl_{12}(\plano)$  of $Cr_1 =0$, $Cr_2= 0$,$ Cr_3 =0$ (considered with multiplicity three), while the  other nine $(-3)$-components are the srict transforms of $P_3 =0$, $Q_3=0$ and $R_3 =0$.

Section \ref{examples}   provide examples  and illustrations of the pencils and of the associated non-minimal elliptic fibration.

\section{Effects  of the  quadratic  Cremona  maps  on the arrangement an on  the pencils $\E_{t}$}

One of the reasons we had to change the notation on the quadratic Cremona maps 
\[\phi_1 = Q_1,\quad  \phi_2= Q_\tau, \quad  \phi_3 = Q_{\tau^2}    \]
was the need to  detail their  effects on the triple points of the dual Hesse arrangement.

Computing their effect on lines of the arrangement, it is not difficult to conclude that, for
\[\mathcal{P}_{12} =   \mathcal{P}_3 (1)  \cup \mathcal{P}_3 (\tau)   \cup \mathcal{P}_3  (\tau^2) \cup \mathcal{P}_ 3 (\infty)\]
as in the previous Section \ref{hesse}, we can verify that
\[\begin{cases} Q_1 (\mathcal{P}_3 (1) )= \mathcal{P}_3 (1), \quad Q_1 (\mathcal{P}_3 (\tau) )= \mathcal{P}_3 (\tau^2), \quad  Q_1 (\mathcal{P}_3 (\infty) )= \mathcal{P}_3 (\infty )\\
Q_\tau (\mathcal{P}_3 (1) )= \mathcal{P}_3 (\tau^2), \quad  Q_\tau (\mathcal{P}_3 (\tau) )= \mathcal{P}_3 (\tau),\quad  Q_\tau (\mathcal{P}_3 (\infty) )= \mathcal{P}_3 (\infty )\\
Q_{\tau^2} (\mathcal{P}_3 (1) )= \mathcal{P}_3 (\tau), \quad Q_{\tau^2} (\mathcal{P}_3 (\tau^2) )= \mathcal{P}_3 (\tau^2), \quad Q_{\tau^2} (\mathcal{P}_3 (\infty) )= \mathcal{P}_3 (\infty )
\end{cases}\]

In general, the strict transform by a quadratic Cremona map  $Q$ with indeterminations $p_1,p_2,p_3$  of a curve $C$ of degree $\mbox{deg}(C)$  is a curve $\overline{C}$ of degree \[2 \mbox{deg}(C) - m(p_{j}, C) - m(p_{k},C) - m(p_{l},C)\]
  with points of multiplicites \[m(p_{j}, \overline{C}) = \mbox{deg}(C) - m_{p_k} - m_{p_l,C), C)} \]   

In the case of $Q_{i}$, for $i= 1,\tau,\tau^2$ we can represent  its effect on curves as follows (we call attention to the  shifts in positions):

\begin{pro}\label{Qisimbolico} 

Let  $C$ be a curve of degree $d$,  having  the same multiplicity $m_1, m_\tau, m_{\tau^2}, m_{\infty} $ at the three points of the  $\mathcal{P}_3 (1), \mathcal{P}_ 3 (\tau), \mathcal{P}_ 3 (\tau^2)$,  respectively. These data represented by the list  \[[ d, m_1 , m_\tau, m_{\tau^2}]\]
Using the same convention,  the strict  transform $\overline{C}_i$  of  $C$ by each quadratic map $Q_i$  has degree and multiplicities:
\[\begin{cases} [ 2 d - 3 m_1,  d -2  m_1,  m_{\tau^2}, m_ \tau],\quad   \mbox{if} \, \,  i= 1\\
[ 2 d - 3 m_\tau,   m_{\tau^2}, d - 2  m_\tau, m_1],\quad   \mbox{if} \,\,  i= \tau\\
[ 2 d - 3 m_{\tau^2},   m_\tau ,  m_1 , d - 2  m_{\tau^2}] ,\quad   \mbox{if}\,\,     i =    \tau^2\\
\end{cases}\]

\end{pro}  

The \emph{strict  transform of a  foliation} on a complex surface  by a birational map is defined by  the pullback of the $1$-form defining the foliation and elimination of curves of singularities. For the birational classification of foliations and the role of elliptic pencils in it we refer  \cite{BGF}.
 
As remarked in \cite{MPEffective}, $Q_1 Q_\tau , Q_{\tau^2}$  preserve the \emph{family}  Lins Neto's foliations $\F_t$, not the individual elements. The main reason for leaving the notations $\phi_1, \phi_2, \phi_3$ and adopting $Q_1, Q_\tau, Q_{\tau^2}$ is the next fact:

\begin{pro}\label{efeitonoparametro}

 The  strict transform of   each  Lins Neto's  foliation $\F_t: \, \Omega + t \cdot \Xi =0$ 
by each quadratic map  $Q_i$, $i \in \{1,\tau,\tau^2\}$ is another element $\F_{q_i(t)}:\quad    \Omega + q_i(t)   \cdot \Xi = 0$, 
where
\[\begin{cases} q_1  (t) = - t -1 \\
q_\tau(t)= -t - \tau\\
q_{\tau^2}(t) = -t - \tau^2\end{cases}\] 
 
\end{pro}

In particular,  these changes in parameter apply to the elliptic pencils $\F_{  m+n \tau } =: \E_{m+n \tau}$.  Next Section \ref{examples} gives  examples.

\section{Examples of effects of  $Q_i$ on elliptic pencils and associated elliptic fibrations}\label{examples}

We analyse here  three examples of elliptic pencils and their modifications by the $Q_i$,  which are  first stages of a  diagram of bifurcations of Cremona maps applied to elliptic pencils in Section \ref{hexalelliptic}.

\begin{exe}\label{sextics}

Let us start with the pencil of elliptic sextics  $\E_{-2}$ of \cite{MPEffective}
\[\E_{-2} : \quad c_1 \cdot l_1 m_1 n_1 \cdot  p^3 +  c_2 \cdot l_2 m_ 3 n_2  \cdot  q^3 = 0\]
which has  a third special element $l_3 m_3 n_2 \cdot  r^3 =0$, where 
\[\mathcal{L}_9 = (l_1 \, l_2 \, l_3 \, m_1\,  m_2\,  m_3 \, n_1 \,  n_2 \, n_3 = 0)\]
and 
\[ p := x+ y + z,\quad q:=  x+ \tau^2 y + \tau z,\quad  r:=  x+ \tau y + \tau^2 z  \]
The pencil has    twelve base-points at   $\mathcal{P}_{12}$. The  generic element is smooth at $\mathcal{P}_{12} \setminus \mathcal{P}_ 3(1)$ and has     ordinary $3$-uple points at the three points of  \[\mathcal{P}_3 (1)= \{(1:1:1), (1: \tau: \tau^2), (1:\tau^2,\tau)\}\]

Figure \ref{novof6g6negative}  illustrates the formation of a non-minimal elliptic fibration after elimination of the twelve base-points by blow ups $\sigma$.

Remark on  the three singular fibers $\tilde{I V}$ of the fibration:  the  $(-3)$-components  are strict transforms of the lines of $\mathcal{L}_9$ ;  the  $(-1)$-components are strict transforms of    lines $p =0$,  $q =0$, $r=0$, respectively. The multiplicity three of each  central element of $\tilde{I V}$ results after  the multiplicity three of $p^3=0$, $q^3=0$, $r^3=0$.

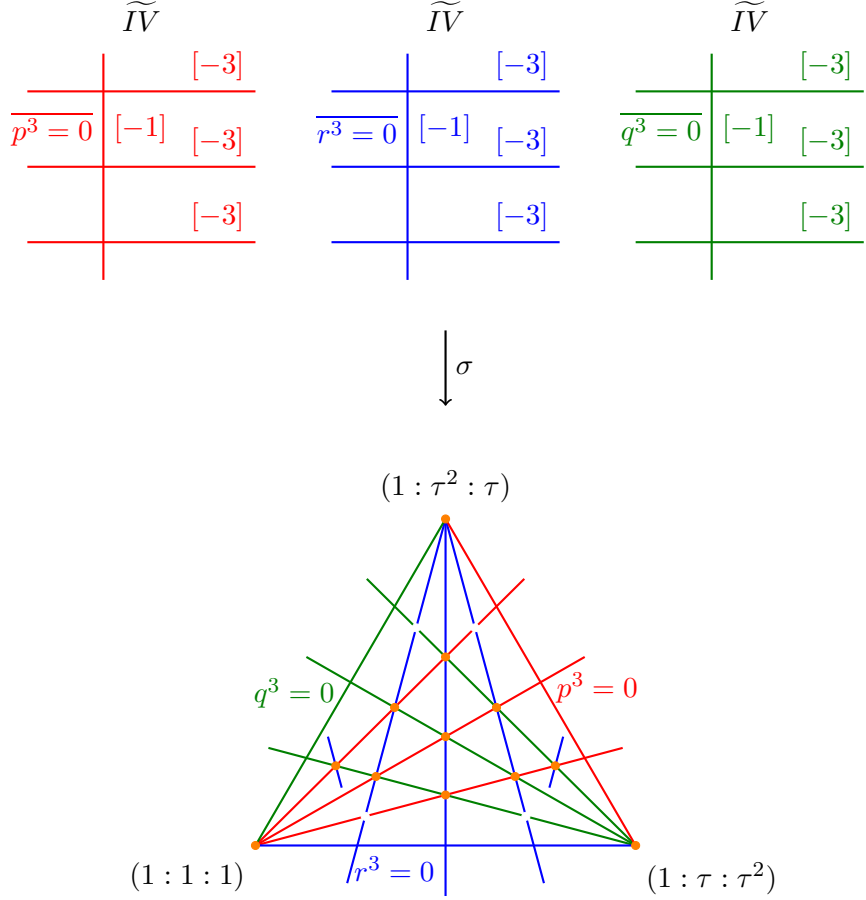
\begin{figure}
   \begin{tikzpicture}[scale=1]
	\coordinate (O) at (0,0);
	\draw[red,thick, name path=O15] (O) -- (15:5);	
	\draw[red,thick, name path=O30] (O) -- (30:5);
	\draw[red,thick, name path=O45] (O) -- (45:5);
	\coordinate (V1) at (5,0);	
	\draw[green!50!black,thick, name path=V115] (V1) -- +(135:5);
	\draw[green!50!black,thick, name path=V130] (V1) -- +(150:5);
	\draw[green!50!black,thick, name path=V145] (V1) -- +(165:5);
	\coordinate (V2) at (canvas polar cs:radius=5cm,angle=60);
	\draw[thick,blue, name path=V215] (V2) -- +(255:5);
	\draw[thick, blue, name path=V230] (V2) -- +(270:5);
	\draw[thick,blue, name path=V245] (V2) -- +(285:5);
	\draw[thick,green!50!black] (O) -- (V2) node[xshift=-2pt,yshift=-4pt,midway,left]{$q^3=0$};
	\draw[thick,blue] (O) -- (V1) node[midway,below left]{$r^3=0$};
	\draw[thick,red] (V1) -- (V2) node[midway,right,xshift=2pt,yshift=-2pt]{$p^3=0$};
	\path [name intersections={of=O15 and V115,by=A1}];
	\path [name intersections={of=O30 and V115,by=A2}];
	\path [name intersections={of=O45 and V115,by=A3}];
	\path [name intersections={of=O15 and V215,by=B1}];
	\path [name intersections={of=O30 and V215,by=B2}];
	\path [name intersections={of=O45 and V215,by=B3}];
	\path [name intersections={of=O30 and V130,by=C1}];
	\path [name intersections={of=O15 and V145,by=C2}];
	\path [name intersections={of=O15 and V130,by=C3}];
	\path [name intersections={of=V215 and V115,by=D1}];
	\path [name intersections={of=V245 and V145,by=D2}];
	\path [name intersections={of=V245 and O45,by=D3}];
	\path [name intersections={of=O45 and V145,by=D4}];
	\draw[blue, thick] ([shift={(-105:0.3cm)}]A1.center) -- ([shift={(-105:-0.4cm)}]A1.center);
	\draw[blue, thick] ([shift={(-75:0.3cm)}]D4.center) -- ([shift={(-75:-0.4cm)}]D4.center);
	\draw[fill,orange] (O) node[black,yshift=-2.5pt,below left]{$(1:1:1)$} circle (1.5pt);
	\draw[fill,orange] (V1) node[black,yshift=-2.5pt,below right]{$(1:\tau:\tau^2)$} circle (1.5pt);
	\draw[fill,orange] (V2) node[black,yshift=2.5pt,above]{$(1:\tau^2:\tau)$} circle (1.5pt);
	\draw[fill,orange] (A1) circle (1.5pt);
	\draw[fill,orange] (A2) circle (1.5pt);
	\draw[fill,orange] (A3) circle (1.5pt);
	\draw[white,fill] (B1) circle (1.5pt);
	\draw[fill,orange] (B2) circle (1.5pt);
	\draw[fill,orange] (B3) circle (1.5pt);
	\draw[fill,orange] (C1) circle (1.5pt);
	\draw[fill,orange] (C2) circle (1.5pt);
	\draw[fill,orange] (C3) circle (1.5pt);
	\draw[white,fill] (D1) circle (1.5pt);
	\draw[white,fill] (D2) circle (1.5pt);
	\draw[white,fill] (D3) circle (1.5pt);
	\draw[fill,orange] (D4) circle (1.5pt);
	\coordinate (C) at ([yshift=2cm] V2);
	\node[right] at (C){$\sigma$};
	\draw[thick,->] ([yshift=0.5cm] C) -- ([yshift=-0.5cm] C);
	\begin{scope}[xshift=-2cm,yshift=8cm,color=red]
	\node[black] at (0.5,3) {$\widetilde{IV}$};
	\node[right] at (0,1.5) {$[-1]$};
	\node[left] at (0,1.5) {$\overline{p^3=0}$};
     \draw[thick] (-1,0) -- (2,0) node[above left]{$[-3]$};
     \draw[thick] (-1,1) -- (2,1) node[above left]{$[-3]$};
     \draw[thick] (-1,2) -- (2,2) node[above left]{$[-3]$};
     \draw[thick] (0,-0.5) -- (0,2.5);
	\end{scope}
	\begin{scope}[xshift=2cm,yshift=8cm,color=blue]
	\node[black] at (0.5,3) {$\widetilde{IV}$};
	\node[right] at (0,1.5) {$[-1]$};
	\node[left] at (0,1.5) {$\overline{r^3=0}$};
	\draw[thick] (-1,0) -- (2,0) node[above left]{$[-3]$};
	\draw[thick] (-1,1) -- (2,1) node[above left]{$[-3]$};
	\draw[thick] (-1,2) -- (2,2) node[above left]{$[-3]$};
	\draw[thick] (0,-0.5) -- (0,2.5);
\end{scope}
	\begin{scope}[xshift=6cm,yshift=8cm,color=green!50!black]
	\node[black] at (0.5,3) {$\widetilde{IV}$};
	\node[right] at (0,1.5) {$[-1]$};
	\node[left] at (0,1.5) {$\overline{q^3=0}$};
	\draw[thick] (-1,0) -- (2,0) node[above left]{$[-3]$};
	\draw[thick] (-1,1) -- (2,1) node[above left]{$[-3]$};
	\draw[thick] (-1,2) -- (2,2) node[above left]{$[-3]$};
	\draw[thick] (0,-0.5) -- (0,2.5);
\end{scope}
\end{tikzpicture}
\caption{The special elements of the pencil of sextics $\F_{-2}$ in red, blue, black  and the nonminimal elliptic fibration obtained after twelve blow ups of base points.}
\label{novof6g6negative}
\end{figure}

\pagebreak

\end{exe}

\begin{exe}\label{nonics}

The strict transform of the pencil of elliptic  sextics $\E_{-2}$  (of previous Example \ref{sextics})  by $Q_\tau$ is a pencil of elliptic nonics  
\[\E_ {2-\tau}   =  \E_{-(-2) -\tau}\]

There are three special elements on this pencil of curves; each one is composed by three lines  of the dual Hesse arrangement  and a triple conic, image by  $Q_\tau$ of a triple line of $\E_{-2}$.

For reference, the conics are 
\[ \begin{cases} 
C_{2,p}= Q_\tau(p)= x^2-\tau xy+y^2+(\tau+1) xz+(\tau +1) yz-(\tau+1) z^2 \\
C_{2,q} =Q_\tau(q)= x^2-xy+\tau y^2-xz+(\tau+1)yz+\tau z^2  \\
C_{2,r} = Q_\tau (r)= x^2+(\tau+1) xy-(\tau+1) y^2-\tau xz+(\tau+1) yz+z^2
\end{cases}\]

The generic element of $\E_ {2-\tau}$  has

\begin{itemize}
\item   three  $4$-uple ordinary points at $\mathcal{P}_3 (\tau) =  \{ (1:1:\tau ), (1: \tau :1), (\tau :1:1)\}$, 
\item three $3$-uple points at $\mathcal{P}_3( \tau^2) =  \{ (\tau^2:1:1), (1:\tau^2:1), (1:1:\tau^2) \}$
(recall that $Q_\tau$ interchanges the sets $\mathcal{P}_3 (1)$  and $\mathcal{P}_3 (\tau^2)$)
\item   smooth points at  the extra six points of $\mathcal{P}_{12}$. 
\end{itemize}

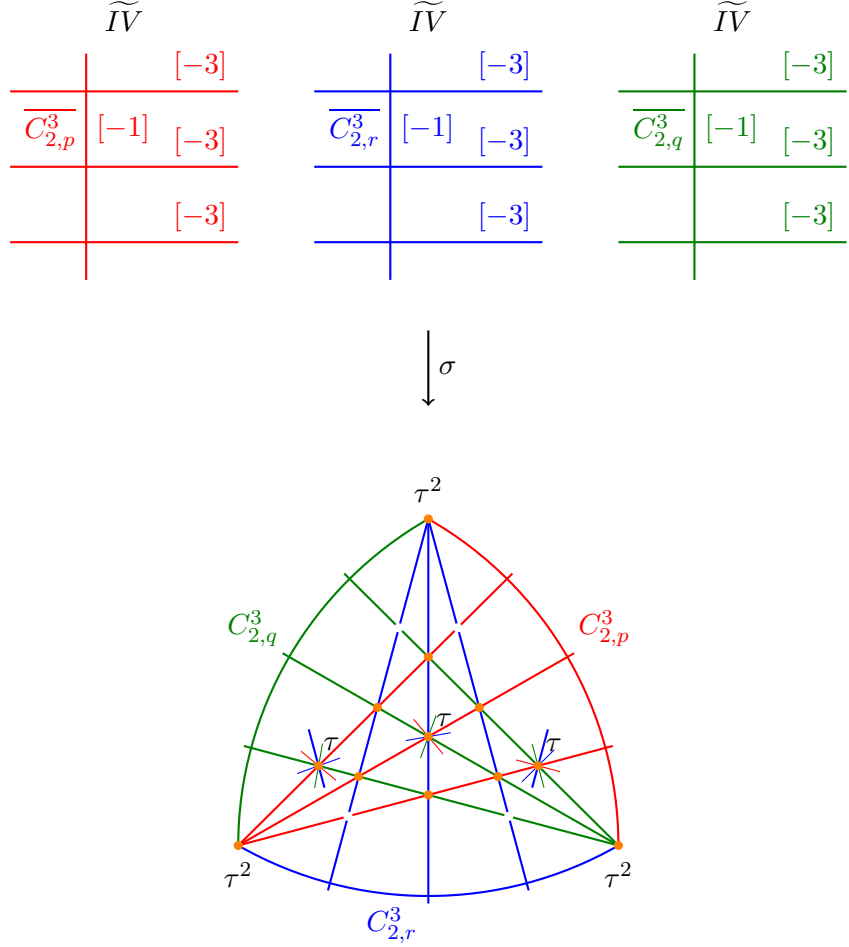
\begin{figure}
   \begin{tikzpicture}[scale=1]
	\coordinate (O) at (0,0);
\draw[red,thick, name path=O15] (O) -- (15:5.1);	
\draw[red,thick, name path=O30] (O) -- (30:5.1);
\draw[red,thick, name path=O45] (O) -- (45:5.1);
\coordinate (V1) at (5,0);	
\draw[green!50!black,thick, name path=V115] (V1) -- +(135:5.1);
\draw[green!50!black,thick, name path=V130] (V1) -- +(150:5.1);
\draw[green!50!black,thick, name path=V145] (V1) -- +(165:5.1);
\coordinate (V2) at (canvas polar cs:radius=5cm,angle=60);
\draw[thick,blue, name path=V215] (V2) -- +(255:5.1);
\draw[thick, blue, name path=V230] (V2) -- +(270:5.1);
\draw[thick,blue, name path=V245] (V2) -- +(285:5.1);
\draw[thick,red] ([shift=(0:5cm)]O) arc (0:60:5cm) node[midway,above right]{$C_{2,p}^3$}; 
\draw[thick,green!50!black] ([shift=(120:5cm)]V1) arc (120:180:5cm) node[midway,above left]{$C_{2,q}^3$};
\draw[thick,blue] ([shift=(240:5cm)]V2) arc (240:300:5cm) node[midway,below left]{$C_{2,r}^3$};
\path [name intersections={of=O15 and V115,by=A1}];
\path [name intersections={of=O30 and V115,by=A2}];
\path [name intersections={of=O45 and V115,by=A3}];
\path [name intersections={of=O15 and V215,by=B1}];
\path [name intersections={of=O30 and V215,by=B2}];
\path [name intersections={of=O45 and V215,by=B3}];
\path [name intersections={of=O30 and V130,by=C1}]; 
\path [name intersections={of=O15 and V145,by=C2}];
\path [name intersections={of=O15 and V130,by=C3}];
\path [name intersections={of=V215 and V115,by=D1}];
\path [name intersections={of=V245 and V145,by=D2}];
\path [name intersections={of=V245 and O45,by=D3}];
\path [name intersections={of=O45 and V145,by=D4}];
\draw[blue, thick] ([shift={(-105:0.3cm)}]A1.center) -- ([shift={(-105:-0.5cm)}]A1.center);
\draw[blue, thick] ([shift={(-75:0.3cm)}]D4.center) -- ([shift={(-75:-0.5cm)}]D4.center);
\foreach \x/\t in {C1/10,A1/45,D4/20} {
	\draw[blue] ([shift={(\t:-0.3cm)}]\x) -- ([shift={(\t:0.3cm)}]\x);
	\draw[green!50!black] ([shift={(\t+60:-0.3cm)}]\x) -- ([shift={(\t+60:0.3cm)}]\x);
	\draw[red] ([shift={(\t+120:-0.3cm)}]\x) -- ([shift={(\t+120:0.3cm)}]\x);
}
\node[above right,xshift=-1pt,yshift=1.5pt] at (A1) {$\tau$};
\node[above right,xshift=-1pt] at (C1) {$\tau$};
\node[above right,xshift=-2pt,yshift=2pt] at (D4) {$\tau$};
\draw[fill,orange] (O) node[black,yshift=-2.5pt,below]{$\tau^2$} circle (1.5pt);
\draw[fill,orange] (V1) node[black,yshift=-2.5pt,below]{$\tau^2$} circle (1.5pt);
\draw[fill,orange] (V2) node[black,yshift=2.5pt,above]{$\tau^2$} circle (1.5pt);
\draw[fill,orange] (A1) circle (1.5pt);
\draw[fill,orange] (A2) circle (1.5pt);
\draw[fill,orange] (A3) circle (1.5pt);
\draw[white,fill] (B1) circle (1.5pt);
\draw[fill,orange] (B2) circle (1.5pt);
\draw[fill,orange] (B3) circle (1.5pt);
\draw[fill,orange] (C1) circle (1.5pt);
\draw[fill,orange] (C2) circle (1.5pt);
\draw[fill,orange] (C3) circle (1.5pt);
\draw[white,fill] (D1) circle (1.5pt);
\draw[white,fill] (D2) circle (1.5pt);
\draw[white,fill] (D3) circle (1.5pt);
\draw[fill,orange] (D4) circle (1.5pt);
	\coordinate (C) at ([yshift=2cm] V2);
	\node[right] at (C){$\sigma$};
	\draw[thick,->] ([yshift=0.5cm] C) -- ([yshift=-0.5cm] C);
	\begin{scope}[xshift=-2cm,yshift=8cm,color=red]
		\node[black] at (0.5,3) {$\widetilde{IV}$};
		\node[right] at (0,1.5) {$[-1]$};
		\node[left] at (0,1.5) {$\overline{C_{2,p}^3}$};
		\draw[thick] (-1,0) -- (2,0) node[above left]{$[-3]$};
		\draw[thick] (-1,1) -- (2,1) node[above left]{$[-3]$};
		\draw[thick] (-1,2) -- (2,2) node[above left]{$[-3]$};
		\draw[thick] (0,-0.5) -- (0,2.5);
	\end{scope}
	\begin{scope}[xshift=2cm,yshift=8cm,color=blue]
		\node[black] at (0.5,3) {$\widetilde{IV}$};
		\node[right] at (0,1.5) {$[-1]$};
		\node[left] at (0,1.5) {$\overline{C_{2,r}^3}$};
		\draw[thick] (-1,0) -- (2,0) node[above left]{$[-3]$};
		\draw[thick] (-1,1) -- (2,1) node[above left]{$[-3]$};
		\draw[thick] (-1,2) -- (2,2) node[above left]{$[-3]$};
		\draw[thick] (0,-0.5) -- (0,2.5);
	\end{scope}
	\begin{scope}[xshift=6cm,yshift=8cm,color=green!50!black]
		\node[black] at (0.5,3) {$\widetilde{IV}$};
		\node[right] at (0,1.5) {$[-1]$};
		\node[left] at (0,1.5) {$\overline{C_{2,q}^3}$};
		\draw[thick] (-1,0) -- (2,0) node[above left]{$[-3]$};
		\draw[thick] (-1,1) -- (2,1) node[above left]{$[-3]$};
		\draw[thick] (-1,2) -- (2,2) node[above left]{$[-3]$};
		\draw[thick] (0,-0.5) -- (0,2.5);
	\end{scope}
\end{tikzpicture}
\caption{At bottom, the special elements of the pencil of nonics $\E_{2-\tau}$ in red, blue, green. Curved lines are representation of conics, represented as not connected. Points denoted by $i$ means points of $\mathcal{P}_3(i)$ At top, three special fibers of type  $\tilde{I V}$ }
\label{novof9g9negative}
\end{figure}

All points of $\mathcal{P}_{12}$ are base-points of the pencil  of nonics The ellimination of base points produces  a non-minimal elliptic fibration, having three special fibers which we shall denote $\tilde{I V}$, having the structure of the fiber $IV$ blown up at its triple points. 

Illustrated by Figure \ref{novof9g9negative}.

Remark on  the $\tilde{I V}$ fibers :  the  $(-3)$-components results from four blow ups at  each projective lines of $\mathcal{L}_9$;  the $(-1)$-components  results from the fact that  the conic components of  $\E_{2-\tau}$ are  blown up five times, namely,  at two  among the three points of $\mathcal{P}_3(\tau^2)$  and at the  three points of $\mathcal{P}_3 (\tau)$.
The multiplicity three of the central elements of $\tilde{I V}$ results after the multiplicity three of the  conic components of $\E_{2-\tau}$.

\end{exe}

\begin{exe}\label{fifteen}

The pencil $\E_{-3+\tau}$ is the strict transfom of $\E_{2-\tau}$ by $Q_1$:
\[\E_{-3+\tau} = \E_{-(2-\tau) -1}\] 

It is a pencil  of elliptic curves of degree fifteen,  
\[\E_{-3+\tau}:\quad  c_1 \cdot P_3 \cdot   (C_{4,p})^3 + c_2\cdot  Q_3  \cdot   (C_{4,q})^3 =0,\]
having a third special element  $R_3 \cdot   (C_{4,r})^3 =0 $, where 
\[\mathcal{L}_9 = ( P_ 3 \cdot Q_3 \cdot  R_3 = 0) \]
and  the  irreducible  rational  quartics  $C_{4,p}, C_{4,q} , C_{4,p}$ are the strict transform by $Q_1$   of the conics $C_p, C_q,C_r$ respectively.

Each  having quartic has   three ordinary $2$-uple points at  $\mathcal{P}_3(1)$. The  generic elliptic  element is a degree fifteen curve  having

\begin{itemize}
\item   three $7$-uple points in $\mathcal{E}_3 (1)$,
\item  three $4$-uple points in $\mathcal{P}_3 (\tau^2)$  (recall that $Q_1$ interchanges $\mathcal{P}_3 (\tau) $ for $\mathcal{P}_3 (\tau^2)$).
\item  three $3$-uple points in $\mathcal{P}_3(\tau)$ 
\item and smooth points at $\mathcal{P}_3 (\infty)$
\end{itemize}


After blow ups  at $\mathcal{P}_{12}$,  the non-minimal fibration has three special fibers $\tilde{IV}$. 
The  $(-3)$-components are  strict transforms of lines of $\mathcal{L}_{12}$. The  $(-1)$-components of the $\tilde{I V}$ are strict transforms of the quartics $C_{4,p}, C_{4,q}, C_{4,r}$;  in fact, each received 3 blow ups at double points and more five  blow ups at smooth points.

\end{exe}

\section{Bifurcations of  Cremona  maps on   elliptic pencils}\label{hexalelliptic} 

\begin{figure}
			\centerline{
			\scalebox{0.8}{
				\xymatrix@C=-0.5cm{ 
					&  &  & \underbrace{[6,3,1,1,1],\textcolor{blue}{1}}_{t=-2} \ar[ld]_-{Q_{\tau}}\ar[rd]^-{Q_{\tau^2}} &  &  &  \\
					&  & \underbrace{[9,1,4,3,1],\textcolor{blue}{2}}_{t=2-\tau}\ar[d]_-{Q_1} &  & \underbrace{[9,1,3,4,1],\textcolor{blue}{2}}_{t=3+\tau}\ar[d]^-{Q_1} &  &  \\
					&  & \underbrace{[15,7,3,4,1],\textcolor{blue}{4}}_{t=-3+\tau} \ar[ld]_-{Q_{\tau}}\ar[rd]^-{Q_{\tau^2}}&  & \underbrace{[15,7,4,3,1],\textcolor{blue}{4}}_{t=-4-\tau}\ar[ld]_-{Q_{\tau}}\ar[rd]^-{Q_{\tau^2}} &  &  \\
					& \underbrace{[21,4,9,7,1],\textcolor{blue}{6}}_{t=3-2\tau}\ar[d]_-{Q_1} &  & \underbrace{[18,3,7,7,1],\textcolor{blue}{5}}_{t=4}\ar[d]^-{Q_1} &  & \underbrace{[21,4,7,9,1],\textcolor{blue}{6}}_{t=5+2\tau}\ar[d]^-{Q_1} &  \\
					& \underbrace{[30,13,7,9,1],\textcolor{blue}{9}}_{t=-4+2\tau}\ar[ld]_-{Q_{\tau}}\ar[rd]^-{Q_{\tau^2}} &  & \underbrace{[27,12,7,7,1],\textcolor{blue}{8}}_{t=-5}\ar[ld]_-{Q_{\tau}}\ar[rd]^-{Q_{\tau^2}} &  & \underbrace{[30,13,9,7,1],\textcolor{blue}{9}}_{t=-6-2\tau}\ar[ld]_-{Q_{\tau}}\ar[rd]^-{Q_{\tau^2}} &  \\
\underbrace{[39,9,16,13,1],\textcolor{blue}{12}}_{t=4-3\tau}\ar[d]_-{Q_1} &  & \underbrace{[33,7,13,12,1],\textcolor{blue}{10}}_{t=5-\tau}\ar[d]_-{Q_1} &  & \underbrace{[33,7,12,13,1],\textcolor{blue}{10}}_{t=6+\tau}\ar[d]^-{Q_1} &  & \underbrace{[39,9,13,16,1],\textcolor{blue}{12}}_{t=7+3\tau}\ar[d]^-{Q_1} \\
\underbrace{[51,21,13,16,1],\textcolor{blue}{16}}_{t=-5+3\tau} &  & \underbrace{[45,19,12,13,1],\textcolor{blue}{14}}_{t=-6+\tau} &  & \underbrace{[45,19,13,12,1],\textcolor{blue}{14}}_{t=-7-\tau} &  & \underbrace{[51,21,16,13,1],\textcolor{blue}{16}}_{t=-8-3\tau}
			}}
		}
\caption{First stages of a diagram of bifurcations of quadratic Cremona maps $Q_i$ applied to elliptic pencils. The general diagram continues the same scheme of applications of the three  $Q_i$ }
\label{diagramelliptic}
\end{figure}
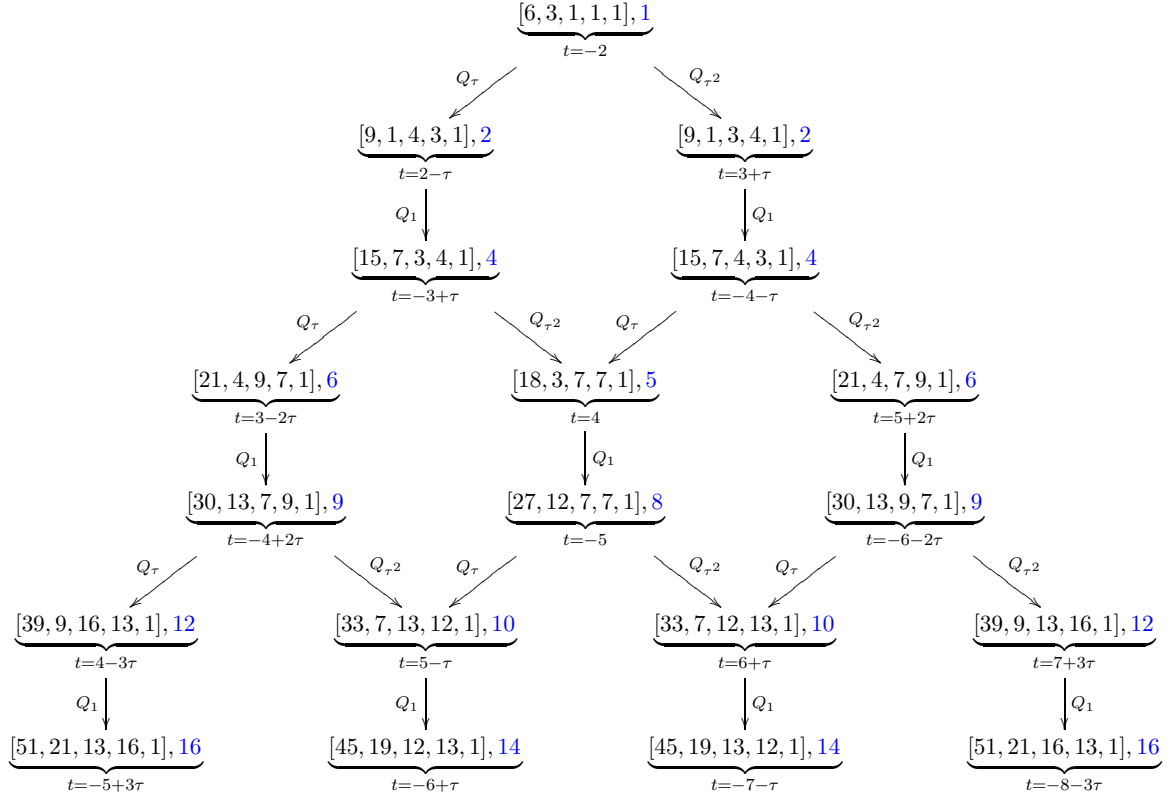

Now we present in Figure \ref{diagramelliptic} the first stages of a  diagram of bifurcations of Cremona maps $Q_1, Q_\tau, Q_\tau^2$ which starts with the  pencil of elliptic sextics $\E_{-2}$  (cf. Example \ref{sextics}).

The  generic element of $\E_{-2}$  is an  elliptic sextic with  three ordinary $3$-uple points at the three points of  \[\mathcal{P}_3 (1)= \{(1:1:1), (1: \tau: \tau^2), (1:\tau^2,\tau)\}\] 
and smooth points at the base points located at  $\mathcal{P}_{12} \setminus \mathcal{P}_ 3(1)$.

We should encode these data as
\[\E_{-2}: \quad [6, \underbrace{[3,3,3]}_{1}, \underbrace{[1,1,1]}_{\tau}, \underbrace{[1,1,1]}_{\tau^2},\underbrace{[1,1,1]}_{\infty}]\]
where $i$ refers to  the set $\mathcal{P}_3 (i)$, $i=1,\tau,\tau^2,\infty$.  For a short notation we can 
use just
\[\E_{-2}: \quad [6, \underbrace{3}_{1}, \underbrace{1}_{\tau}, \underbrace{1}_{\tau^2},\underbrace{1}_{\infty}]\] or even simpler 
\[ \E_{-2}: \quad [6,3,1,1,1],    \]
keeping the order $i, \tau,\tau^2,\infty$, But when using the  formula for  geometrical genus  or when  determining the self-intersection after blow ups at base-points,  it is important to  remember that each multiplicity represents multiplicities at  three points. 

So in  general, 
 \[[d, m_1, 
m_{\tau}, m_{\tau^2}, m_\infty]\]
means the degree $d$  of the generic elliptic curve and $m_i$ means multiplicity of singulairites \emph{at the three points} of the sets $\mathcal{P}_3 (i)$ , for $i= 1,\tau,\tau^2, \infty$,  in this order.

\newpage 

For obtaning  concretely  any data \[[d, m_1, 
m_{\tau}, m_{\tau^2}, m_\infty]\]   of  the diagram, for any row and column, use the next   code in \emph{Python}. In the default example, the code   determines the data $[45,19,12,13,1]$ of  entry $a_{7,2}$, which is:
\[ Q_1(Q_ {\tau^2}(Q_1(Q_\tau (Q_1(Q_\tau([6,3,1,1,1]))))))\]
by printing 
\[  \mbox{q1(q3(q1(q2(q1(q2(list))))))}\]

\begin{lstlisting} 
list=[6,3,1,1,1] # qinf([3,1,1,1,0])
def q1(list):  #q1
    l0=list[0];l1= list[1];l2=list[2];l3=list[3];
    list[0]=2*l0-3*l1
    list[1]= l0-2*l1
    list[2]=l3
    list[3]=l2
    return(list)
def q2(list): # q_tau
    l0=list[0];l1= list[1];l2=list[2];l3=list[3];
    list[0]= 2*l0-3*l2
    list[1]=l3
    list[3]=l1
    list[2]=l0-2*l2
    return(list)
def  q3(list): #q_tau^2
    l0=list[0];l1= list[1];l2=list[2];l3=list[3];
    list[0]= 2*l0-3*l3
    list[1]= l2
    list[2]=l1
    list[3]=l0-2*l3
    return(list)
print(q1(q3(q1(q2(q1(q2(list)))))))

\end{lstlisting}

The extra  blue number at right is the  degree of   each one of the three  rational components of $\E_{t} $ having multiplicity three.    For instance, 
\[ \E_{-2}: \quad [6,3,1,1,1],  \color{blue}  1 \] 
means one straight line  having multiplicity three as special component of $\E_{-2}$.

The degree of  each rational plane curve $Cr$, component of multiplicity three  of the elliptic pencil
\[[d,m_1,m_\tau,m_{\tau^2},m_\infty]\] 
 in the diagram Figure \ref{diagramelliptic}  is
\[\mbox{deg}(C_r) = \frac{d}{3} -1\]
In  fact, as already remarked,   each pencil in the diagram  is of the  
\[\E_{m+n\tau}: \quad   c_1 \cdot P_3 \cdot  (Cr_1)^3 + c_2 \cdot  Q_3 \cdot  (Cr_2)^3  = 0\]
where $P_3=0$ and  $Q_3=0$ define sets of three distinct straight lines belonging to $\mathcal{L}_9$.  Hence  $d = 3 + 3 \mbox{deg}(Cr_{j}) $.

Now we begin describing  the  structure of 
the infinite  diagram whose first stages are in Figure \ref{diagramelliptic}.

We shall use repeatedly the relation
\[ 1 + \tau + \tau^2 = 0\]

\begin{lem}\label{t(i,j)}

The  parameter $t(i,j)$ of the elliptic pencil $\E_{t(i,j)}$ encoded in the entry    $a_{i,j}$ ($i, j\in \mathbb{N}$) in the infinite diagram of bifurcations whose first stages are given in Figure \ref{diagramelliptic} verify
\[t(1,1)   = -2,\quad   t(2,1)= 2- \tau ,\quad  t(3,1) = - 3 +\tau;\]
for all $i \geq 2$ and $j=1$: 
\[\begin{cases}   t(2i,1)  = t(2,1)  +  (i-1) (1-\tau) \\
  t(2i+1,1)  = t(3,1)  -  (i-1) (1-\tau) \end{cases}\]
for all $i \geq 2$ and $2\leq  j\leq i+1$:
  \[\begin{cases}  t(2i,j)= t(2i,1) + (j-1)(1+2\tau) \\  t(2i+1,j)= t(2i+1,1) - (j-1)(1+2\tau) \end{cases}\]

\end{lem}

Proof:

Let us suppose the formula for  $j=1$ and  up to some  $i\geq 2$. Then $t(2 (i+1),1 )= t(2i +2, 1)$ 
is obtained after the composition $q_\tau \circ q_1$  applied to $t(2i, 1)$, that is,
\[ t(2 (i+1),1 )= q_\tau( q_1 (t(2i ,1)))= \]
\[=  - (- t(2i ,1) -1) -\tau =\]
\[= t(2i ,1 ) + 1 -\tau =\]
\[= t(2,1) + (i-1) (1-\tau) + 1 - \tau =\]
\[=  t(2,1) + i (1-\tau)\]
which is the formula for $i+1$. 

Also for $j=1$ and up to $i\geq 2$, $t(2 (i+1) + 1,1 )= t(2i +1 + 2 , 1)$ is obtained after the composition $q_1 \circ q_\tau $  applied to $t(2i +1, 1)$, that is,   
\[ t(2 (i+1) + 1,1 )= q_1 ( q_\tau(t(2i +1, 1) )) =\]
\[= - (- t(2i +1,1) -\tau) -1 =\]
\[= t(2i +1,1 ) + \tau -1 =\]
\[= t(3,1) - (i-1) (1-\tau) - ( 1 - \tau) =\]
\[=  t(3,1) - i (1-\tau)\]
which is the formula for $i+1$. 

Now supose the formula for $i$ and up to $j\geq 2$. 

Remark that, in the scheme of the diagram,   $t(2 i, j+1)$ is obtained after applying to $t(2 i , j)$ the composition
 \[ q_{\tau^2} \circ q_{\tau}^{-1} =   q_{\tau^2} \circ q_\tau, \]
hence
\[t(2 i , j+1) =
 - (- t(2 i, j) -\tau) - \tau^2 =\]
\[= t(2 i, j) + \tau -\tau^2=\]
\[= t(2 i, j) + 1 + 2 \tau    =\] 
\[=  t(2 i ,1) + (j-1) (1+2\tau) + 1+ 2\tau=\]
\[=  t(2 i ,1) + j (1 + 2\tau)\]
 which is the formula for $j+1$. 
 
Again we suppose the formula up to $i\geq 2$. Remark that $t(2 i + 1 , j+1)$ is obtained after applying to $t(2 i +1 , j)$ the composition  \[q_1 \circ  q_{\tau^2} \circ q_{\tau}^{-1} \circ q_1^{-1} =  q_1 \circ  q_{\tau^2} \circ q_\tau \circ q_1, \]
hence
\[t(2 i + 1, j+1) = - (-( - (- t(2 i +1, j) -1) - \tau ) - \tau^2)- 1  =\]
\[= t(2 i + 1 , j) -\tau + \tau^2 =\]
\[= t(2 i + 1 , j) -1 - 2\tau =\]
\[= t(2 i + 1 , 1) - (j-1) (1+2\tau) - ( 1+ 2\tau)=\]
\[=  t(2 i +1  ,1) - j (1 + 2\tau)\]
 which is the formula for $j+1$.  
 
\qed

\begin{lem}\label{cong1}

All parameters $t(i,j) = m+ n \tau$   described in Lemma \ref
{t(i,j)}   verify
\[m+n \equiv 1 \, \mbox{mod}(3)\]

\end{lem}

Proof:

For the initial parameter   $t(2,1) = 2 -\tau$, it holds  $2-1 \equiv 1 \, \mbox{mod}(3) $;   for the initial parameter   $t(3,1) = -3  + \tau$ also it  holds  $-3 + 1  \equiv 1 \, \mbox{mod}(3) $.

According to Lemma \ref{t(i,j)},   throughout the diagram are added  to the initial parameters   terms  of the form 
\[ m^{\prime} + n^{\prime} \tau := \pm (i-1) (1-\tau)\]   or terms of the form
\[  m^{\prime} + n^{\prime} \tau := \pm (j-1) (1+ 2 \tau) \]
In both cases, 
\[m^{\prime} +n^{\prime} \equiv  0 \,\,  \mbox{mod}(3)\]  
\qed

The fact indicated in Lemma \ref{cong1}  is important for the correct application  of the formula of \cite{Puchuri} on the degrees of the generic element of the elliptic pencil $\E_{t(i,j)}$. 

In fact, $t(i,j)= m+ n \tau$ with $m +n \equiv 1 \, \mbox{mod} (3) $ is   equivalent to 
\[\mbox{gcd}( t(i,j)- 1, -2 -\tau) \neq 1\]
The main Theorem of \cite{Puchuri} demands  the simplification of the quotient 
\[ \frac{t(i,j)- 1}{ -2 -\tau} \in \mathbb{Z}(\tau),\]
theme of the next Lemma \ref{simplification}.

\begin{lem}\label{simplification}

For any parameter $t(i,j)$ as in Lemma \ref{t(i,j)}   it holds:
\[\frac{t(i,j) - 1}{-2-\tau} =   c(i,j) + d(i,j) \tau \in \mathbb{Z} (\tau),\]
namely,   $c(1,1)=1, d(1,1) = - 1$;

if  $ i\geq  2$  and  $i$ is even:
\[\begin{cases} \mbox{if}\,\,     j=1: \quad  c(i,j) =0, d(i,j)= \frac{i}{2}\\
   \mbox{if} \, \,   2 \leq    j \leq  \frac{i}{2}  +1: \quad  c(i,j)= - (j-1) , d(i,j)= \frac{i}{2} + 1- j \end{cases}\]

if $ i \geq 2$ and  $i$ is odd:  
\[\begin{cases} \mbox{if} \, \, j= 1: c(i,j)=1,  d(i,j)=  - \frac{i-1}{2}   -1\\ 
  \mbox{if} \, \,  2 \leq    j \leq \frac{i-1}{2}  +1:\quad  c(i,j) = 1 +  (j-1) ,  d(i,j) = -\frac{i-1}{2}  -1 + (j-1) \end{cases}\] 
\end{lem}

Proof:

We start  with
\[\frac{t(1,1)- 1}{-2-\tau} = \frac{-3}{-2-\tau} = 1- \tau =: c(1,1) + d(1,1) \tau\]
Supose $j=1$ and for \emph{even} $i$ and  up to some  $i\geq 2$ the validity of: 
\[ \frac{t(i ,1) -1}{-2-\tau} = \frac{i}{2} \tau\]
It follows from  Lemma \ref{t(i,j)} that
\[t(i + 2, 1) =  q_\tau (q_1 (t(i, 1)))  =\]
\[=  t(i ,1) + (1-\tau),\]
hence
\[ \frac{ t(i + 2, 1) -1 }{-2-\tau} = \frac{t(i ,1) + (1-\tau) - 1}{-2-\tau} =\]
\[= \frac{t(i,1)-1}{-2-\tau} + \frac{1- \tau}{-2-\tau} =\] 
\[ = \frac{i}{2} \tau + \tau = \frac{i+2}{2} \tau,\]
which is the formula for the next even $i+2$. 

Supose $j=1$ and for \emph{odd} $i$  up to  $i\geq 2$  the validity of:
\[ \frac{t(i ,1) -1}{-2-\tau} =  1 - (\frac{i-1}{2} + 1) \tau\]
It follows from  Lemma \ref{t(i,j)} that
\[t(i + 2, 1) = t(i ,1) - (1-\tau),\]
hence
\[ \frac{ t(i + 2, 1) -1 }{-2-\tau} = \frac{t(i ,1) - (1-\tau) - 1}{-2-\tau} =\]
\[= \frac{t(i,1)-1}{-2-\tau} - \frac{1-\tau}{-2-\tau}  =\] 
\[ = 1 - (\frac{i-1}{2} + 1) \tau - \frac{1- \tau}{-2-\tau}  =\]
\[= 1 - (\frac{i-1}{2} + 1) \tau - \tau =  \]
\[= 1  - (\frac{i+1}{2} + 1) \tau  \]
which is the formula for the next odd $i+2$. 

Suppose the formula for even $i$  and up to $j\geq 2$. it follows from  Lemma \ref{t(i,j)}  that
\[ t(i,j+1) = t(i,j) + (1+2\tau)\]Hence
\[\frac{t(i,j+1) -1}{-2-\tau} = \frac{t(i,j)-1}{-2-\tau} + \frac{1+2 \tau}{-2-\tau} =\]
\[=  - (j-1) + (\frac{i}{2} + 1 -j) \tau +   \frac{1+2 \tau}{-2-\tau} =  \]
\[=  - (j-1) + (\frac{i}{2} + 1 -j) \tau    - (1+\tau) =   \]
\[= - j + (\frac{i}{2} + 1 - (j+1) ) \tau,\]
which is the formula for $j+1$.

Suppose the formula for odd $i$  and up to $j\geq 2$. it follows from  Lemma \ref{t(i,j)}  that
\[ t(i,j+1) = t(i,j) - (1+2\tau)\]Hence
\[\frac{t(i,j+1) -1}{-2-\tau} = \frac{t(i,j)-1}{-2-\tau} - \frac{1+2 \tau}{-2-\tau} =\]
\[=  1 +  (j-1) +  (-\frac{i-1}{2} - 1 + (j-1)) \tau -   \frac{1+2 \tau}{-2-\tau} =  \]
\[= 1 +  (j-1) +  (-\frac{i-1}{2} - 1 + (j-1)) \tau  + (1+\tau) =   \]
\[= 1 +  j + (-\frac{i-1}{2} - 1 +j) \tau,\]
which is the formula for $j+1$.

\qed

\begin{pro}\label{cdt}

The degree of the rational curve $Cr_{i,j}$ encoded in entry $a_{i,j}$  ($i, j\in \mathbb{N}$) in the infinite diagram of bifurcations whose first stages are given in Figure \ref{diagramelliptic} verify
\[ c^2(i,j) + d^2(i,j) -c(i,j) + d(i,j) - c(i,j) d(i,j)  \]
where $c(i,j)$ and $d(i,j)$ are given in Lemma \ref{simplification}. 
\end{pro}

Proof:

After Lemma \ref{cong1} and Lemma \ref{simplification} we are in conditions to apply the formula of \cite{Puchuri}(Theorem, p. 142)  for the degree of the generic element of the elliptic pencil $\E_{t(i,j)}$ as   
{\small \[  d_{t(i,j)}  = 3 \cdot (\, a^2+b^2+c^2(i,j)+d^2(i,j)- a b-a c(i,j) + a d(i,j)  - b d(i,j) - c(i,j) d(i,j) \, ) \]}
with  $a=1$ and $b=0$ and $c(i,j)$, $d(i,j)$ given in Lemma \ref{simplification}.

As already remarked  \[\mbox{deg}(Cr_{i,j})= \frac{d_{t(i,j)}}{3} -1\] 

After simplification,  we get the result. 
\qed 

Therefore we have proved Theorem \ref{fatoprincipal}. 

The result is implemented as the next  \emph{Python} code. The  default choice $a_ {11,6}$ is the last entry of the sequence of blue numbers given below, associated to degree $36$. 

\begin{lstlisting}

aij=[13,7]
i= aij[0];j=aij[1];
a=1;b=0;
def deg(a,b,c,d):
    print((a**2+b**2+c**2+d**2-a*b-a*c+a*d-b*d-c*d )-1)
if i>=1:
    if  i==1 and  j==1:
        c=1; d=-1;
        deg(a,b,c,d);
    elif i >= 2:
        if   (i % 2) == 0:
            if j <=0 or   j > i/2 +1:
                print('For a[i,j] with even i,  choose j with 1 <= j <= i/2 +1')
            else:
                if  j==1: 
                    c=0 ; d= i/2;
                    deg(a,b,c,d);
                elif 2<=   j and j  <= i/2 +1:  
                    c= - (j-1) ;  d= i/2 + 1- j ;
                    deg(a,b,c,d);
        elif  (i % 2) ==1:
            if j <=0 or   j > (i-1)/2 +1:
                print('For a[i,j] with odd i, choose j with 1 <= j <= (i-1)/2 + 1')
            if   j== 1: 
                c=1; d= - (i-1)/2  -1;
                deg(a,b,c,d);
            elif  2<=   j and j  <= (i-1)/2 +1:
                c= 1+  (j-1) ;  d= -(i-1)/2 -1 + (j-1);
                deg(a,b,c,d);
else:
     print('Choose  a[i,j] with  i>=1')                    







\end{lstlisting}

\begin{center}
\begin{tabular}{>{$}l<{$\hspace{12pt}}*{13}{c}}
&&&&&&&\color{blue} 1&&&&&&\\
&&&&&&\color{blue} 2&&\color{blue} 2&&&&&\\
&&&&&&\color{blue} 4&&\color{blue} 4&&&&&\\
&&&&& \color{blue} 6 &&\color{blue} 5&& \color{blue} 6 &&&&\\
&&&&& \color{blue} 9 &&\color{blue} 8&& \color{blue} 9 &&&&\\
&&&& \color{blue} 12 &&\color{blue} 10&& \color{blue} 10 && \color{blue} 12  &&&\\
&&&& \color{blue} 16 &&\color{blue} 14&& \color{blue} 14 && \color{blue} 16  &&&\\
&&& \color{blue} 20 &&\color{blue} 17&& \color{blue} 16 && \color{blue} 17 &&\color{blue} 20 &&\\
&&& \color{blue} 25 &&\color{blue} 22&& \color{blue} 21 && \color{blue} 22 &&\color{blue} 25 &&\\
&& \color{blue} 30 &&\color{blue} 26&& \color{blue} 24 && \color{blue} 24 &&\color{blue} 26  &&\color{blue} 30 &\\
&& \color{blue} 36 &&\color{blue} 32&& \color{blue} 30 && \color{blue} 30 &&\color{blue} 32  &&\color{blue} 36 &\\
\end{tabular}
\end{center}
\centerline{\color{blue} etc.}

\section{Singularities and equations of of rational curves in the diagram}\label{rationalcomponentshexalelliptic}

In \emph{Table  1}  of \cite{preprint} the authors start listing degrees and singularities of  rational curves appearing in their bifurcations of  Cremona maps.

They list: the straight line;  a smooth conic; a  quartic with three $2$-uple points; a   sextic with four $2$-uple points and two $3$-uple points; a  nonic  with four $2$-uple, two $3$-uple and three $4$-uple points;  a curve of  degree twelve with four $2$-uple, two $3$-uple, three $4$-uple and three $5$-uple points.

Now we can give a complete  description of any  rational curves appearing in our Cremona bifurcations, describing  degree (blue) and   multiplicities of singular points at  each set of   points $P_3(i)$, $i=1,\tau,\tau^2$

The data   
	\[ [{\color{blue} d},[m_{1,1}, m_{1,2}, m_{1,3}],[m_{\tau,1}, m_{\tau,2}, m_{\tau,3}],[m_{\tau^2,1}, m_{\tau^2,2}, m_{\tau^2,3}] \]
  give in  blue the degree of the  the rational curve  and each \emph{non-ordered} $3$-tuple give the multiplicities in the three points of the sets $P_3(i)$, $i=1,\tau,\tau^2$.
  
In our choice, the starting straight line and all of its strict transforms  by any sequence of applications of by $Q_1, Q_\tau,Q_\tau^2$ \emph{do not pass} by points of the set $\mathcal{P}_3 (\infty)$. 

For example,  
\[ [{\color{blue} 4} , [2, 2, 2], [1,1,1], [1,1,0]]
 \]
means an irreducible rational quartic, having  three $2$-uple points at $P_3(1)$,  passing simply by the three points of   $P_3(\tau)$,  and passing simply by two of the three points of   $P_3(\tau^2)$ (no order considered)

\begin{figure}
	\centerline{
	\scalebox{0.6}{
		\xymatrix@C=-1cm{
			&  &  & [\textcolor{blue}{1},[1,1,0],[0,0,0],[0,0,0]] \ar[ld]_{Q_{\tau}}\ar[rd]^{Q_{\tau^2}} &  &  &  \\
			&  & [\textcolor{blue}{2},[0,0,0],[1,1,1],[1,1,0]]\ar[d]_{Q_1} &  & [\textcolor{blue}{2},[0,0,0],[1,1,0],[1,1,1]]\ar[d]^{Q_1} &  &  \\
			&  & [\textcolor{blue}{4},[2,2,2],[1,1,0],[1,1,1]] \ar[ld]_{Q_{\tau}}\ar[rd]^{Q_{\tau^2}}&  & [\textcolor{blue}{4},[2,2,2],[1,1,1],[1,1,0]]\ar[ld]_{Q_{\tau}}\ar[rd]^{Q_{\tau^2}} &  &  \\
			& [\textcolor{blue}{6},[1,1,1],[3,3,2],[2,2,2]]\ar[d]_{Q_1} &  & [\textcolor{blue}{5},[1,1,0],[2,2,2],[2,2,2]]\ar[d]^{Q_1} &  & [\textcolor{blue}{6},[1,1,1],[2,2,2],[3,3,2]]\ar[d]^{Q_1} &  \\
			& [\textcolor{blue}{9},[4,4,4],[2,2,2],[3,3,2]]\ar[ld]_{Q_{\tau}}\ar[rd]^{Q_{\tau^2}} &  & [\textcolor{blue}{8},[4,4,3],[2,2,2],[2,2,2]]\ar[ld]_{Q_{\tau}}\ar[rd]^{Q_{\tau^2}} &  & [\textcolor{blue}{9},[4,4,4],[3,3,2],[2,2,2]]\ar[ld]_{Q_{\tau}}\ar[rd]^{Q_{\tau^2}} &  \\
			[\textcolor{blue}{12},[3,3,2],[5,5,5],[4,4,4]]\ar[d]_{Q_1} &  & [\textcolor{blue}{10},[2,2,2],[4,4,4],[4,4,3]]\ar[d]_{Q_1} &  & [\textcolor{blue}{10},[2,2,2],[4,4,3],[4,4,4]]\ar[d]^{Q_1} &  & [\textcolor{blue}{12},[3,3,2],[4,4,4],[5,5,5]]\ar[d]^{Q_1} \\
			[\textcolor{blue}{16},[7,7,6],[4,4,4],[5,5,5]] &  & [\textcolor{blue}{14},[6,6,6],[4,4,3],[4,4,4]] &  & [\textcolor{blue}{14},[6,6,6],[4,4,4],[4,4,3]] &  & [\textcolor{blue}{16},[7,7,6],[5,5,5],[4,4,4]]
	}}
}
\caption{Bifurcations of quadratic Cremona maps $Q_i$ applied to irreducible  rational curves. }
\label{diagramrational}
\end{figure}
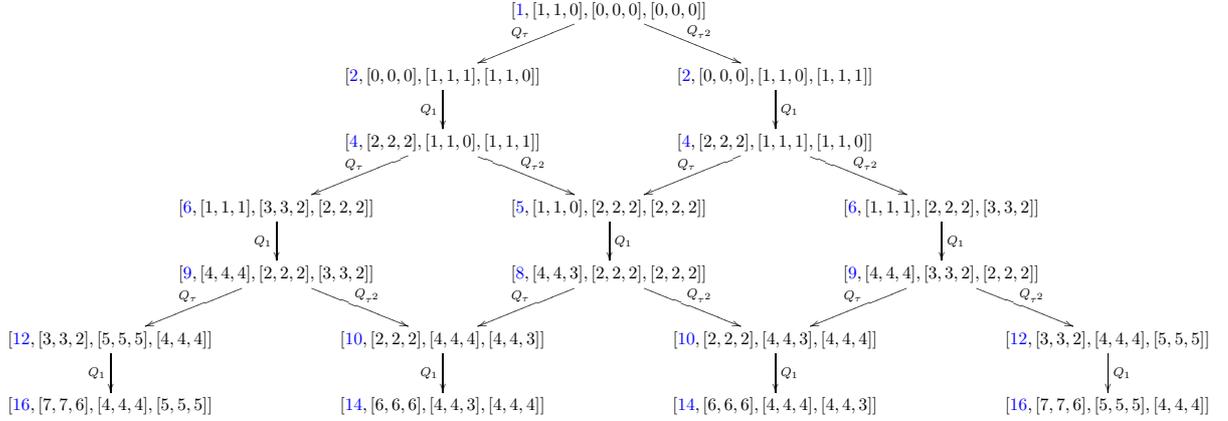

\pagebreak

The data of degree and singularities at each  position in the diagram can be computed by the next code in \emph{Python}.   The entry $a_{6,3}$ 
\[ [{\color{blue} 10}, [2, 2, 2], [4, 4, 3], [4, 4, 4]]\]
is the default example in the next code:

\begin{lstlisting}
list=[1,[1,1,0],[0,0,0],[0,0,0]] 
def q1(list):  #q1
    l0=list[0];l10= list[1][0];l11=list[1][1]; l12=list[1][2]; l2=list[2];l3=list[3];
    list[0]=2*l0-l10-l11-l12
    list[1]= [l0-l11-l12, l0-l10-l12, l0-l10-l11]
    list[2]=l3
    list[3]=l2
    return(list)
def q2(list): # q_tau
    l0=list[0];l1= list[1];l20= list[2][0];l21=list[2][1]; l22=list[2][2]; l3=list[3];
    list[0]= 2*l0-l20-l21-l22
    list[1]=l3
    list[3]=l1
    list[2]=[l0-l21-l22, l0-l20-l22, l0-l20-l21]
    return(list)
def  q3(list): #q_tau^2
    l0=list[0];l1= list[1];l2=list[2];l30= list[3][0];l31=list[3][1]; l32=list[3][2];
    list[0]= 2*l0-l30-l31-l32
    list[1]= l2
    list[2]=l1
    list[3]=[l0-l31-l32, l0-l30-l32, l0-l30-l31]
    return(list)
print(q3(q1(q3(q1(q2(list))))))



\end{lstlisting}

Explicit equations of the rational curves can be obtained easily with the \emph{Singular} software. for instance the quintic at entry $a_{4,2}$ is the strict transform of the line $x+y+z=0$ by the composition $Q_{\tau^2} \circ Q_1 \circ Q_\tau$ and has equation

{\small  \[Cr_{4,2}:  x^5-x^4y-2x^3y^2-2x^2y^3-xy^4+y^5-x^4z-x^3yz-3 x^2y^2z-xy^3z-y^4z-2x^3z^2-\]
\[- 3x^2yz^2-3xy^2z^2-2y^3z^2-2x^2z^3-xyz^3-2y^2z^3-xz^4-yz^4+z^5   = 0\]}
This is the default example in the next code for the  \emph{Singular} software, which produces a factorization of the total transform.    

\begin{lstlisting}
ring R=(0,a),(x,y,z),dp;                                      
minpoly=a2+a+1;                                                     
poly f = x+y+z;                                                                                     
map Q1 = R, y^2- x*z, x^2- y*z, z^2- x*y;                                                        
map Qt = R, a*y^2- x*z, a*x^2- y*z, z^2- a^2 *x*y;                                               
map Qt2 = R, a^2*y^2- x*z, a^2 *x^2- y*z, z^2- a*x*y; 
map cre= Qt;
cre=Q1(cre);
cre=Qt2(cre);
poly cref=cre(f);                                          
factorize(cref);
\end{lstlisting}

\section{Acknowledgements}

We are very grateful to Dr. Orestes Bueno for his help on the  design  of  figures and diagrams.

The second author was partially supported by \textsc{pucp-Peru}~(\textsc{dgi}: 2023-E-0020).

{Liliana Puchuri, Pontificia Universidad Cat\'olica del Per\'u, Per\'u.

email: lpuchuri@pucp.pe

Lu\'{i}s  Gustavo Mendes, Universidade Federal do Rio Grande do Sul, UFRGS, Brazil. 

email:   gustavo.mendes@ufrgs.br}

\end{document}